\documentclass[12pt]{article}
\usepackage[cp1251]{inputenc}
\usepackage[russian, english]{babel}
\usepackage[psamsfonts]{amssymb}
\usepackage{amsfonts,euscript}
\usepackage{amsmath,amsthm,graphicx}
\usepackage{color}
\usepackage[colorlinks]{hyperref}
\usepackage{pgf,tikz,pgfkeys}
\usepackage{mathrsfs}
\usetikzlibrary{arrows.meta,calc}

\newcommand{\midarr}{\tikz \draw[-{Latex[length=4mm, width=2mm]}] (0,0) -- +(.3,0);}
\newcommand{\midrra}{\tikz \draw[-{Latex[length=4mm, width=2mm]}] (0,0) -- +(-.3,0);}

\newtheorem{thm}{\bf Theorem}
\newtheorem{lem}[thm]{\bf Lemma}

\newcommand{\rr}{\mathbb{R}}

\newcommand {\al} {\alpha}

\newcommand {\da} {\delta}
\newcommand {\Da} {\Delta}

\newcommand {\sa} {\sigma}

\newcommand {\fy} {\varphi}

\newcommand{\IN}{{\subset}}

\newcommand {\mmm}{{\setminus}}

\newcommand{\8}{{\infty}}
\newcommand{\io}{{I^\infty}}
\newcommand{\ia}{{I^*}}
\newcommand{\0}{{\varnothing}}
\newcommand{\vse}{$\blacksquare$}

\newcommand{\bj}{{\bf {j}}}

\newcommand{\eS}{{\EuScript S}}

\newcommand{\eC}{{\EuScript C}}
\newcommand{\eP}{{\EuScript P}}

\newcommand{\eK}{{\EuScript K}}
\newcommand{\eF}{{\EuScript F}}

\def \Lip {\mathop{\rm Lip}\nolimits}
\def \min {\mathop{\rm min}\nolimits}

\title{Even unique intersection point can break OSC: an example.}

\author{Kirill Kamalutdinov \and 
Andrey Tetenov}

\begin{document}

\maketitle


\begin{abstract} 
This was a long-standing question since 90-ies whether one-point intersection property for a self-similar set implies open set condition. We answer this question negatively.
We give an example of a totally disconnected self-similar set  $K\IN\rr$ which does not have   open set condition  and has minimal overlap of its pieces, that is, all intersections of its pieces $K_i\cap K_j,i\neq j$ are empty except only one, which is a single point. 
\end{abstract}


\subsection{Introduction}

The aim of current short paper is to prove the following:
\begin{thm}
There exist a system $\eS=\{S_1,S_2,..,S_6\}$ of contraction similarities in $\rr$ with totally disconnected attractor $K\IN[0,1]$, $\{0,1\}\IN K$  for which the following holds:\\
1. $\eS$ does not satisfy OSC;\\
2. The only non-empty intersection of the pieces of $K$ is $K_3\cap K_4= S_3(0)=S_4(1)$ so the system $\eS$ is post-critically finite;\\
3. The Hausdorff dimension $\dim_H(K)$ is equal to similarity dimension $s$ of the system $\eS$, while its measure $H^s(K)=0$.\\
 \end{thm}
 
 It is well-known \cite{Hut} that if a system $\eS=\{S_1,...,S_m\}$ of contraction maps in complete metric space $X$ has the attractor $K$ for which all $K_i=S_i(K)$ are disjoint, it satisfies open set condition (OSC). Nevertheless, it seemed unclear how OSC relates to the actual size of the overlap  \cite{SSS7,BR}. Since early 90-ies there were many attempts to prove that one point intersection property implies OSC. It was proved by C.Bandt and H.Rao in \cite{BR} that a finite overlap implies OSC
for connected self-similar sets in the plane, 
 while  our example  \cite{TKV} shows that there are families of  self-similar arcs in $\rr^3$,  satisfying one-point intersection property, which do not satisfy OSC.

\newpage

 So the question arises what is the situation in disconnected case. It was shown in \cite{BR}, that there are Cantor sets of arbitrary small dimensions which do not fulfill the OSC. However, in this case one of the overlapping points should have a recurrent address, so these sets are not post-critically finite.

There are several works \cite{PolSim,Bar} which produce the families of Cantor sets with overlaps, but neither of these examples guarantee that the intersections of their pieces are finite. These works are based on the use of potential-theoretic characterization of Hausdorff dimension, and on the transversality condition, which   first appeared in \cite{PolSim, Solom}.

In the current paper we construct a system 
$\eS=\{S_1,..,S_4\}$ of contraction similarities in $\rr$ with the attractor $K$, for which all the pieces $K_i$ are disjoint except $K_2$ and $K_3$, whose intersection is a singleton,  whereas  OSC, and even weak separation property  fail to hold for the system $\eS$. Moreover, both addresses of the overlapping point are periodic, the set $K$ is of finite type and is post-critically finite.

The cornerstone for this construction is our approach, based on General Position Theorem \cite[Theorem 14]{TF}. It allows us to construct  families of self-similar sets with prescribed behavior of their critical sets.   It gave exact overlap for double fixed points in \cite{TF}, while in \cite{TKV} it allowed to obtain one-point intersections for the pieces of self-similar arc in $\rr^3$ which does not satisfy OSC. In the current work this method is applied to prove the existence of systems in $\rr$ with unique one point intersection, not satisfying WSP.

There is one more aspect of the presented example. The Hausdorff dimension of the attractor $K$ of the system $\eS$ in the Theorem 1 is equal to the similarity dimension $s$ of the system $\eS$.  At the same time, by \cite{Schief,Zer}, the set $K$ should have zero Hausdorff measure in dimension $s$. The peculiarity of the situation is that the intersections $K_i\cap K_j$ of different pieces  of the attractor is at most one point, therefore the measure drop cannot be caused by the overlap of the pieces, but only by their relative position.

The authors are grateful to Vladislav Aseev, Christoph Bandt and Caroly Simon, whose comments were of great value for the authors.

\subsection{Some preliminaries.}

Let $\eS=\{S_1,\dots,S_m\}$ be a system of contraction similarities in $\mathbb R^n$. A nonempty compact set $K=K(\eS)$ such that $K=\bigcup\limits_{i=1}^{m} S_i (K)$, is called an {\em attractor} of the system $\eS$, or a {\em self-similar set} generated by the system $\eS$.

By $I=\{1,2,...,m\}$ we denote the set of indices, $\ia=\bigcup\limits_{n=1}^\8 I^n$  
is the set of all  multiindices $\bj=j_1j_2...j_n$. So  $F=\{S_\bj, \bj\in\ia\}$ will denote the semigroup, generated by $\eS$. 
The set of all infinite sequences $I^{\8}=\{{\bf \al}=\al_1\al_2\ldots,\ \ \al_i\in I\}$ is  the
{\em index space}; and $\pi:I^{\8}\rightarrow K$ is the {\em index map
 }, which sends a sequence $\bf\al$ to  the point $\bigcap\limits_{n=1}^\8 K_{\al_1\ldots\al_n}$.

The system $\eS$ is said to satisfy the {\em open set condition} (OSC), if there exists an open set $O$ such that $S_i(O)\subset O$ and $S_i(O)\cap S_j(O)=\varnothing$ for all distinct $i,j\in I=\{1,\dots,m\}$.

Denote by $F=\{S_{\bf i}: {\bf i}\in I^\infty\}$ the semigroup, generated by $\eS$; then $\eF=F^{-1}\circ F$, or a set of all compositions $S_{\bf j}^{-1}S_{\bf i}$, ${\bf i}, {\bf j} \in I^{*}$, is the {\em associated family of similarities}. The system $\eS$ has the {\em weak separation property} (WSP) iff $\rm Id \notin \overline{\eF \setminus \rm Id}$. If the system doesn't have WSP, then it doesn't satisfy OSC, but the opposite is not true.

 A set $C(\eS)=\bigcup\limits_{i=1,j\neq i}^{m} S_i (K)\cap S_j (K)$ is called a {\em critical set} of the system $\eS$.
 $\eS$ is called  { postcritically finite} or PCF, if the set 
$\eP = \{ \al \in \io : \exists  i_1 \ldots i_n : S_{i_1} \ldots  S_{i_n} ( \pi ( \al ) ) \in \eC \} $
is finite.

\subsection{The construction.}

Take $p,q,r$ in $(0,  1/36)$ and put  $h =\dfrac{8}{15},a=\dfrac{3}{15}$.
Define a system $\eS_{pqr}=\{S_1,S_2,...,S_6\}$ of contraction similarities of $[0,1]$ depending on parameters $p,q,r$ by the equations 
$$S_1 (x)=p x,\quad S_2 (x)=a +rx,\quad S_3 (x)=h -qx, \quad S_4 (x)=h-r+rx,$$ $$ \quad S_5 (x)=1-a -rx,\quad S_6 (x)=1-r+rx ,$$ \\
Let  $K_{pqr}$ be the attractor of the system $\eS_{pqr}$ and $K_i=S_i(K_{pqr})$ be its pieces. Notice that $K_{pqr}\IN[0,1]$ and $\{0,1\}\IN K_{pqr}$.
By the construction, the only possible non-empty intersection of the pieces of $K_{pqr}$ is $K_3\cap K_4$ which always contains the point $h$.

For simplicity, we write $K$ and $\eS$ instead of $K_{pqr}$ and $\eS_{pqr}$ if it does not cause ambiguity.\\

\vspace{1cm}

\hspace{-1.3cm}\begin{tikzpicture}[line cap=round,line join=round,>=stealth,x=1.1cm,y=1.1cm]
\draw[line width=.75 pt,dotted] (-3,0.1)--(12,0.1);
\draw [line width=1.pt] (-3.,0)--node[blue] {\midarr}  (-1.6,0);
\draw [line width=2pt,green] (0.4,0)--node[blue] {\midarr}  (1.7,0);
\draw [line width=2pt,green] (8.7,0)--node[blue] {\midrra}  (7.4,0);
\draw [line width=1.5pt] (10.5,0)--node[black] {\midarr}(12.,0);

\draw [line width=1.5pt] (5.3,0.2)--node[red] {\midrra}(4.3,0.2);
\draw [line width=1.5 pt] (3.8,0)--node[black] {\midarr}(5.3,0);
\draw [fill=black] (-3,0.05) circle (2.2pt)(5.3,0.1) circle (3pt)(12,0.05) circle (2.2pt);
\foreach \position in {(-1.6,0.), (0.4,0),(1.7,0),(7.4,0), (8.7,0), (10.5,0.), (3.8,0), (4.3,0.2)} \draw[fill=black] \position circle (1.5pt);
\draw (-2.5,0.35) node {$K_1$}(1,0.35) node {$K_2$} (8,0.35) node {$K_5$}(11.5,0.35) node {$K_6$};
\draw (-2.5,-0.25) node {$px$}(1,-0.25) node {$rx$}(4.7,0.45) node {$K_3; -qx$} (4.5,-0.25) node {$K_4;rx$}(8,-0.25) node {$-rx$}(11.5,-0.25) node {$rx$} (5.4,-0.18) node {$h$}(-3.2,-.2) node {\large\bf 0} (12.3,0) node {\large\bf 1};

\node [below=.7cm, align=flush center,text width=8cm] at (4,0)
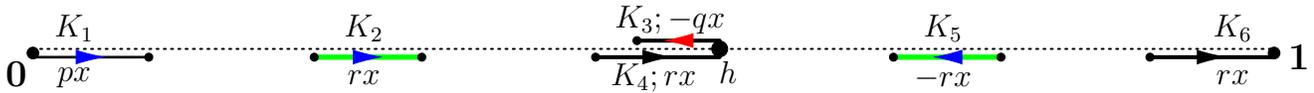
        {\small Relative position of the pieces of $K$.
            
        };
\end{tikzpicture}

\vspace{1cm}

Our aim is to find such values of $p,q,r$ that $S_2(K)\cap S_3 (K)=\{h\}$, and we will say in this case that the system $\eS$ has {\em unique one-point intersection}.\\

{\bf Remark.}  The system $\eS_{pqr}$ may have both unique one point intersection property and OSC. For example, it happens for any $p=r\in(0,1/36)$ and $\dfrac{144}{175}<q<\dfrac{7}{8}$ because  in this case $S_3(K\mmm K_1)\cap S_4(K\mmm K_6)=\0$.\\

\begin{lem}\label{nowsp} If $\dfrac{\log p}{\log r}\notin\mathbb Q$, then the system $\eS_{pqr}$  does not have WSP for any $q$.
\end{lem}
{\bf Proof:  }
Denote $H_m(x)=S_3 S_1^m S_5(x)$ and $G_n(x)=S_4 S_6^n S_2(x)$.
Routine computation gives  $H_m(x)=h-p^mq +p^mqa-p^mqrx$ and
  and 
$G_n(x)=h-r^{n+1}+r^{n+1}a-r^{n+2}x$.\\

It is a well-known fact then, (cf. \cite[Lemma 7]{TF}) that for any $q>0$ there is a sequence  $(m_k,n_k)\in \mathbb N^2$, such that $ p^{-m_k}r^{n_k+1}$ converges to  $q$ as $k\to\infty$.

  Consider the sequence $$G_{n_k}^{-1} H_{m_k} (x) = \dfrac{(r^{n_k+1}-p^{m_k}q)(1-a)}{r^{n_k+2}} + \dfrac{p^{m_k}q}{r^{n_k+1}} x.$$  Since $\lim\limits_{k\to\8}\dfrac{r^{n_k+1}-p^{m_k}q}{r^{n_k+1}}=0$ and
$\lim\limits_{k\to\8}\dfrac{p^{m_k}q}{r^{n_k+1}}=1$, the sequence $G_k^{-1} H_k$ converges to identity. So the point $\rm Id$ is a limit point of the associated family $ \mathcal{F}$, which contradicts WSP. \vse \ \\

\subsection{\bf General position and displacement theorems}

To prove the existence of the parameters $p,q,r$ for which $\eS$ has unique one-point intersection, we use a connective of two statements, which we called General Position Theorem and Displacement Theorem. We refer the reader to \cite{TF} or \cite{TKV} for the proofs of these theorems.
First is General Position Theorem, which we use in less general form than in \cite{TF}:\\

\smallskip

\begin{thm}\label{genpos} 
Let $(D,\rho),(L_1,\sa_1), (L_2,\sa_2)$  be   metric spaces. \\
Let \ \  $\varphi_i(\xi,x):\ D\times L_i\to \mathbb R^n$  be  continuous maps,  such that\\
 {\bf (a)} they  are $\alpha$-H\"older  with  respect  to $x$;
 and\\  
 {\bf (b)} there is  $M>0$ such that  for  any $x_1\in L_1, x_2\in L_2$, $\xi,\xi'\in D$\\ the  function \quad  $\Phi(\xi,x_1,x_2)=\varphi_1(\xi,x_1)-\varphi_2(\xi,x_2)$ \quad satisfies  the inequality 
\begin{equation}\label{antilip}\left\|\Phi(\xi',x_1,x_2)-\Phi(\xi,x_1,x_2)\right\|\ge M\rho(\xi',\xi)    \end{equation}Then Hausdorff dimension of the set $\Delta=\{\xi\in D|\ \   \varphi_1(\xi,L_1)\cap\varphi_2(\xi,L_2)\neq\0\} $ satisfies \begin{equation}\label{dimda }\dim_H \Delta \leq \min\left\{\dfrac{\dim_H L_1\times L_2}{\alpha},d\right\}\end{equation} 
Moreover, if the spaces $L_1,L_2$ are compact, the set $\Da$ is closed in $D$.
\end{thm}

\vspace{.5cm} 
In our situation, the set $D$ will be some interval in $\rr$, $L_1, L_2$ will be the address space $I^\8$, and $\fy_i$ will be the maps of the type $S_\bj\circ\pi$, sending $I^\8$ to parametrized pieces of $K$.\\

To evaluate the displacement  $|\pi (\sigma)-\pi' (\sigma)|$ of elements $x=\pi (\sigma)$ of the set $K_{pqr}$ under the transition to the set $K_{pq'r}$ caused by change of the parameter $q$, we use the following Displacement Theorem:

\medskip

\begin{thm}\label{collage}

Let  $\eS=\{S_1,...,S_m  \}$ and $\eS'=\{S'_1,...,S'_m  \}$   be  two  systems of  contractions in $\rr^n$. Let  $\pi:\ I^\infty\to K$ and  $\pi':\ I^\infty\to K'$ be the address maps with $I=\{1,...,m\}$. Suppose $V$ is  such compact set, that for  any $i=1,...,m$, $S_i(V)\IN V$  and  $ S'_i(V)\IN V$.\\
Then, for  any $\sa\in \io$, 
\begin{equation}\label{deltapi}\|\pi(\sa)-\pi'(\sa)\|\le\dfrac{\da}{1-R},\end{equation}
where   $R=\max\limits_{i\in I} (\Lip S_i,\Lip S'_i)$  and $\da=\max\limits_{x\in V, i\in I} \| S'_i(x)- S_i(x)\|$.
\end{thm}

\subsection{\bf Applying general position theorem}

Now we  wish to evaluate the set of those $p,q,r$, for which $K_3 \cap K_4=\{h\}$. 
First we fix
 $p,r$, and let $q$ to vary so that we could apply  Theorem \ref{genpos} to evaluate the dimension of the set
\begin{equation}\label{baddmn}
\Delta(p,r)=\left\{q\in (0,1/36):\ K_3 \cap K_4 \neq \{h\}\right\}
\end{equation}


Notice that
\begin{equation}\label{union}
K=\{0\}\cup\bigcup\limits_{m=0}^\infty S_1^m (K \setminus K_1)=\{1\}\cup\bigcup\limits_{n=0}^\infty S_4^n (K \setminus K_6),\end{equation}
therefore $K_3$ and $K_4$ can be represented as:
$$K_3=\{h\}\cup\bigcup\limits_{m=0}^\infty S_3 S_1^m (K \setminus K_1),\ \ K_4=\{h\}\cup\bigcup\limits_{n=0}^\infty S_4 S_6^n (K \setminus K_6)$$

Then  $K_3 \cap K_4=\{h\}$  iff for any $m,n\in \mathbb N\cup\{0\}$ and any $i\in I\setminus\{1\}$,  $j\in I\setminus \{6\}$, $S_3 S_1^m (K_i)\cap S_4 S_6^n (K_j)= \varnothing$.\\


To apply Theorem \ref{genpos}, for each $i\in I\setminus\{1\}$ and  $j\in I\setminus \{6\}$ we consider the functions  and   $\varphi_q(\sa)=S_3 S_1^m S_i \pi_{pqr}(\sa)$ and $\psi_q(\tau)=S_4 S_6^n S_j \pi_{pqr}(\tau)$ acting from $I^\8$ to $K$. \\

In view to satisfy the condition (a) of the Theorem \ref{genpos} we supply the space  $I^\infty$ with a metrics in which these functions are Lipshitz:

{\bf  The space $I^\8_R$.} \quad Let $0<R<1$ and $I^\8_R$ be the space $I^\8$ supplied with the metrics $\rho_R(\sa,\tau)=R^{w(\sa,\tau)}$, where $w(\sa,\tau)=\min\{k:\sa_k=\tau_k\}-1$.\\
 This metrics turns $I^\8$ to a self-similar set having Hausdorff dimension $\dim_H I^\8_R=-\dfrac{\log 6}{\log R}$. Particularly, if $0<R<\dfrac{1}{36}$,\ then $\dim_H I^\8_R<1/2$.\\
\begin{lem}\label{tech1}
Let $p,q,r\in(0,R]$, $R\in(0,1)$. Then the map $\pi_{pqr}:\ I^\8_R\to K_{pqr}$ is  1-Lipschitz. \vse \end{lem}

Now we check the condition (b) of Theorem \ref{genpos}.

 Note that if  $i\in I\mmm\{1\}$, $j\in I\mmm\{6\}$, and $S_3 S_1^m S_i (K) \cap S_4 S_6^n S_j (K) \neq \0$, then $S_3 S_1^m [a,1] \cap S_4 S_6^n [0, 1-a] \neq \0$ which is equivalent to  $p^m [aq,q] \cap r^{n+1} [a,1] \neq \0$. Then the inequality
 $r^{n+1} \le p^m \dfrac{q}{a}$ should hold, which means \begin{equation}\label{dmnpr}\dfrac{r^{n+1}}{p^m}\le  \dfrac{q}{a} \end{equation} 

So, in search of those $q$ for which $S_3 S_1^m S_i (K)$ and  $S_4 S_6^n S_j (K)$ may intersect, we can restrict the values of $q$ to the intervals
$$D_{mn}(p,r):=\left(\dfrac{a r^{n+1} }{p^m}, r \right)$$ 

These intervals will serve as the parameter set $D$ in the Theorem \ref{genpos}. Now we need to check that there is an inequality similar to  (\ref{antilip}) for properly chosen maps:
%
%
%

\begin{lem}\label{tech2}
Let $i\in I\mmm\{1\} $, $j\in I\mmm\{6\}$. Let $\varphi_q(\sa)=S_3 S_1^m S_i \pi_{pqr}(\sa)$ and $\psi_q(\tau)=S_4 S_6^n S_j \pi_{pqr}(\tau)$. Then for any $\sa,\tau\in I^\8$ and for any $q,q' \in D_{mn}(p,r)$: 
\begin{equation}\label{ge11d }|\varphi_q(\sa)-\psi_q(\tau)-\varphi_{q'}(\sa)+\psi_{q'}(\tau)|>\dfrac{p^m}{35}  |q-q'| \end{equation}
\end{lem}

{\bf Proof:  }\\

Take  $\eS=\{S_1,...,S_6\}$ and $\eS'=\{S'_1,...,S'_6\}$, 

Let $x=S_i \pi_{pqr}(\sa)$, $x'=S'_i \pi_{pq'r}(\sa)$, $y=S_j \pi_{pqr}(\tau)$, $y'=S'_j \pi_{pq'r}(\tau)$ be the images of $\sa,\tau$ in $K$ and $K'$.\\

Let $F=S_3 S_1^m$, $G=S_4 S_6^n$, $F'=S'_3 S_1^m$, $G'=G$. Denote $\da=|q-q'|$.\\

It follows from Theorem \ref{collage} that $|x-x'|$ and $|y-y'|$ do not exceed $\dfrac{36\da}{35}$. Since $x=S_i \pi_{pqr}(\sa)$, where $i\neq 1$, we have $x\ge a$.

Taking into account the inequality (\ref{dmnpr}), we have:\\

$|F(x)-G(y)-F'(x')+G'(y')|=\left|p^m (q x - q' x) + p^m (q' x - q' x') + r^{n+1} (y - y')\right|$


$\ge p^m \left( |q- q'| x - q' |x - x'| - \dfrac{r^{n+1}}{p^m} |y - y'|\right) > p^m \left( a - \dfrac{1}{36} \cdot \dfrac{36}{35} - \dfrac{q}{a}\cdot\dfrac{36}{35}\right)\da$\\

\bigskip

Finally, using that $a=1/5$ and $\dfrac{q}{a}<\dfrac{5}{36}$, we get:\\

$|F(x)-G(y)-F'(x')+G'(y')|>\dfrac{p^m}{35} \da$
\vse \\

\subsection{\bf Almost all $K_{pqr}$ have unique one point intersection.}



Let $\Delta_{mn}(p,r)$ denote the set of all $q\in D_{mn}(p,r)$ such that for some $i\in I\setminus\{1\}$ and $j\in I\setminus\{6\}$,\qquad
  $S_3 S_1^m S_i (K_{pqr}) \bigcap S_4 S_6^n S_j (K_{pqr})) \neq \0$.

\begin{lem}\label{tech3}
Let $p\in(0,r)$, $r\in(0,1/36)$. Then for any $m,n\in\mathbb N$ the set $\Delta_{mn}(p,r)$ is a closed subset of $D_{mn}(p,r)$ whose dimension is at most $-\dfrac{2\log 6}{\log r}$ .
\end{lem}
{\bf Proof:  }\ 
Take some  $i\in  I\setminus\{1\}$, $j\in I\setminus\{6\}$ and consider the functions $\fy_1(q,\sa)=S_3 S_1^m S_i  \pi_{pqr}(\sa)$ and $\fy_2(q, \sa)= S_4 S_6^n S_j \pi_{pqr}(\sa)$ mapping $I^\8_r$ to $K_{pqr}$. It follows from Lemma \ref{tech1} that both these functions are Lipschitz with respect to $\sa$, and from Lemma \ref{tech2} it follows that if $q,q'\in D_{mn}(p,r)$ and $\Phi(q,\sa,\tau)=\fy_1(q,\sa)-\fy_2(q,\tau)$ then: 
\begin{equation}\label{deltaF}|\Phi(q',\sa,\tau)-\Phi(q,\sa,\tau)|\ge \dfrac{p^m}{35}|q'-q|\end{equation}

Applying Theorem \ref{genpos} to $D_{mn}(p,r)$ we get that the set $$\Da^{(ij)}_{mn}(p,r)=\{q\in D_{mn}(p,r): S_3 S_1^m S_i (K_{pqr}) \bigcap S_4 S_6^n S_j (K_{pqr})) \neq \0\}$$ is closed in $D_{mn}(p,r)$ and its dimension is not greater than $2\dim_H I^\8_r<1$. The set $\Da_{mn}(p,r)$ is the union of all $\Da^{(ij)}_{mn}(p,r)$ so it  is  closed   in $D_{mn}(p,r)$ and has dimension not greater than $2\dim_H I^\8_r<1$.
\vse 

\medskip

Since the set $\Delta(p,r)$ is a countable union
 $\bigcup\limits_{m,n=0}^{\infty} \Da_{mn}(p,r)$ of the sets whose Hausdorff dimension  is not greater than $2\dim_H I^\8_r$, the same is true for $\Delta(p,r)$. 

Therefore $\Delta(p,r)$ has zero Lebesgue measure in $\rr$.

Finally, by Fubini's Theorem, we have:

\begin{thm}\label{fullmes}
For each $r\in(0,1/36)$ the set $\eK_r$ of those ${(p,q)}\in (0,r)^2$, for which $\eS_{pqr}$ has unique one point intersection, has full measure in $(0,r)^2$, and for 
each $p\in (0,r)$  the set $\eK_r$ has full measure in $(0,r)$.
\end{thm}

\subsection{\bf Dimension calculation.}

\begin{thm}\label{hdim} If $\eS_{pqr}$ has unique one point intersection, then Hausdorff dimension $d=\dim_H K_{pqr}$ satisfies the equation \qquad  $p^d + q^d+ 2r^d=1$.
\end{thm}
{\bf Proof:  }\

This is obvious in the case when $\eS_{pqr}$ satisfies OSC.\\

Otherwise, note that the set $K$ may  be considered as the attractor of an infinite system $\eS^{*}=\{S_1^k S_j:\ k\in \mathbb N\cup\{0\},\ j\in\{2,3,4\}\}$. Thus we have that $d=\dim_H K_{pq}$ has upper estimate: $d\le d^{*}$, where $d^{*}$ is  the unique solution of the equation 
\begin{equation}\label{eq}\sum\limits_{k=0}^{\infty} p^{kd} \left(q^d+2r^d\right)=1,\end{equation} 
which is equivalent to $p^d + q^d+2r^d=1$.

Consider the sequence of subsystems $\eS_n = \{S_1^k S_j:\ k\in \{0,1,\dots,n\},\ j\in\{2,3,4\}\}$ in the system $\eS^{*}$. Let $K_n$ be an attractor of $\eS_n$. Obviously $K_n \subset K$, so $d_n=\dim_H K_n \le d$. From the other side, since $0\notin K_n$, we have that $S_2(K_n)\cap S_3 (K_n)=\varnothing$. Therefore the system $\eS_n$ satisfy OSC and the dimension $d_n$ of it's attractor is the unique solution of the equation $\sum\limits_{k=0}^{n} p^{kd} \left(q^d+2r^d\right)=1$.

The sequence $d_n$ increases and $d_n\le d$, so it has a limit which satisfies the equation \ref{eq}, i.~e. $d^{*}=\lim\limits_{n\to\infty}d_n$. Therefore $d=d^{*}$. \quad \vse

Finally, if $\eS_{pqr}$ does not satisfy OSC and has unique one point intersection property, then $H^d(K)$ cannot be positive by Schief's Theorem \cite{Schief}, so $H^d(K)=0$.

\end{document}